

 \magnification=1200 \hsize=11.25cm \vsize=18cm \parskip 0pt \parindent=12pt
 \voffset=1cm \hoffset=1cm


\catcode'32=9

\font\tenpc=cmcsc10
\font\eightpc=cmcsc8
\font\eightrm=cmr8
\font\eighti=cmmi8
\font\eightsy=cmsy8
\font\eightbf=cmbx8
\font\eighttt=cmtt8
\font\eightit=cmti8
\font\eightsl=cmsl8
\font\sixrm=cmr6
\font\sixi=cmmi6
\font\sixsy=cmsy6
\font\sixbf=cmbx6

\skewchar\eighti='177 \skewchar\sixi='177
\skewchar\eightsy='60 \skewchar\sixsy='60

\catcode`@=11

\def\tenpoint{%
  \textfont0=\tenrm \scriptfont0=\sevenrm \scriptscriptfont0=\fiverm
  \def\rm{\fam\z@\tenrm}%
  \textfont1=\teni \scriptfont1=\seveni \scriptscriptfont1=\fivei
  \def\oldstyle{\fam\@ne\teni}%
  \textfont2=\tensy \scriptfont2=\sevensy \scriptscriptfont2=\fivesy
  \textfont\itfam=\tenit
  \def\it{\fam\itfam\tenit}%
  \textfont\slfam=\tensl
  \def\sl{\fam\slfam\tensl}%
  \textfont\bffam=\tenbf \scriptfont\bffam=\sevenbf
  \scriptscriptfont\bffam=\fivebf
  \def\bf{\fam\bffam\tenbf}%
  \textfont\ttfam=\tentt
  \def\tt{\fam\ttfam\tentt}%
  \abovedisplayskip=12pt plus 3pt minus 9pt
  \abovedisplayshortskip=0pt plus 3pt
  \belowdisplayskip=12pt plus 3pt minus 9pt
  \belowdisplayshortskip=7pt plus 3pt minus 4pt
  \smallskipamount=3pt plus 1pt minus 1pt
  \medskipamount=6pt plus 2pt minus 2pt
  \bigskipamount=12pt plus 4pt minus 4pt
  \normalbaselineskip=12pt
  \setbox\strutbox=\hbox{\vrule height8.5pt depth3.5pt width0pt}%
  \let\bigf@ntpc=\tenrm \let\smallf@ntpc=\sevenrm
  \let\petcap=\tenpc
  \normalbaselines\rm}

\def\eightpoint{%
  \textfont0=\eightrm \scriptfont0=\sixrm \scriptscriptfont0=\fiverm
  \def\rm{\fam\z@\eightrm}%
  \textfont1=\eighti \scriptfont1=\sixi \scriptscriptfont1=\fivei
  \def\oldstyle{\fam\@ne\eighti}%
  \textfont2=\eightsy \scriptfont2=\sixsy \scriptscriptfont2=\fivesy
  \textfont\itfam=\eightit
  \def\it{\fam\itfam\eightit}%
  \textfont\slfam=\eightsl
  \def\sl{\fam\slfam\eightsl}%
  \textfont\bffam=\eightbf \scriptfont\bffam=\sixbf
  \scriptscriptfont\bffam=\fivebf
  \def\bf{\fam\bffam\eightbf}%
  \textfont\ttfam=\eighttt
  \def\tt{\fam\ttfam\eighttt}%
  \abovedisplayskip=9pt plus 2pt minus 6pt
  \abovedisplayshortskip=0pt plus 2pt
  \belowdisplayskip=9pt plus 2pt minus 6pt
  \belowdisplayshortskip=5pt plus 2pt minus 3pt
  \smallskipamount=2pt plus 1pt minus 1pt
  \medskipamount=4pt plus 2pt minus 1pt
  \bigskipamount=9pt plus 3pt minus 3pt
  \normalbaselineskip=9pt
  \setbox\strutbox=\hbox{\vrule height7pt depth2pt width0pt}%
  \let\bigf@ntpc=\eightrm \let\smallf@ntpc=\sixrm
  \let\petcap=\eightpc
  \normalbaselines\rm}
\catcode`@=12

\tenpoint
\font\tenbboard=msbm10
\def\bboard#1{\hbox{\tenbboard #1}}

\catcode`\@=11
\def\pc#1#2|{{\bigf@ntpc #1\penalty \@MM\hskip\z@skip\smallf@ntpc%
    \uppercase{#2}}}
\catcode`\@=12

\def\pointir{\discretionary{.}{}{.\kern.35em---\kern.7em}\nobreak
   \hskip 0em plus .3em minus .4em }

\def\qed{\quad\raise -2pt\hbox{\vrule\vbox to 10pt{\hrule width 4pt
   \vfill\hrule}\vrule}}

\def\rem#1|{\par\medskip{{\it #1}\pointir}}

\def\vspace[#1]{\noalign{\vskip#1}}

\def\resume#1{\vbox{\eightpoint\narrower\narrower
\pc R\'ESUM\'E|.\quad #1}}
\def\abstract#1{\vbox{\eightpoint\narrower\narrower
\pc ABSTRACT|.\quad #1}}

\def\section#1{\goodbreak\par\vskip .3cm\centerline{\bf #1}
   \par\nobreak\vskip 3pt }

\long\def\th#1|#2\endth{\par\medbreak
   {\petcap #1\pointir}{\it #2}\par\medbreak}

\def\article#1|#2|#3|#4|#5|#6|#7|
    {{\leftskip=7mm\noindent
     \hangindent=2mm\hangafter=1
     \llap{{\tt [#1]}\hskip.35em}{\petcap#2}\pointir
     #3, {\sl #4}, {\bf #5} ({\oldstyle #6}),
     pp.\nobreak\ #7.\par}}

\def\livre#1|#2|#3|#4|
    {{\leftskip=7mm\noindent
    \hangindent=2mm\hangafter=1
    \llap{{\tt [#1]}\hskip.35em}{\petcap#2}\pointir
    {\sl #3}, #4.\par}}

\def\divers#1|#2|#3|
    {{\leftskip=7mm\noindent
    \hangindent=2mm\hangafter=1
     \llap{{\tt [#1]}\hskip.35em}{\petcap#2}\pointir
     #3.\par}}

\font\sevenbboard=msbm7

\def\setA{{\bboard A}}

\def\setAA{\hbox{\sevenbboard A}}
\def\setZ{{\bboard Z}}

\def\setFBA{{\cal A}}
\def\matB{{\bf B}}
\def\matF{{\bf F}}
\def\matI{{\bf I}}
\def\matM{{\bf M}}

\def\len{{\ell}}
\def\inv{{\rm inv}} 
%
\def\Ferm{\mathop{\rm Ferm}\nolimits}
\def\Bos{\mathop{\rm Bos}\nolimits}
\def\Univ{\mathop{\rm Univ}\nolimits}
\def\bw#1#2{{#1\choose #2}}

\def\smallmatrix#1{\vcenter{\offinterlineskip
     \halign{\vrule height 4pt depth 2pt width 0pt
      \hfil$\scriptstyle##$\hfil&&\kern
       3pt\hfil$\scriptstyle##$\hfil \crcr#1\crcr}}}

\catcode`\@=11
\def\matrice#1{\null \,\vcenter {\normalbaselines \m@th
\ialign {\hfil $##$\hfil &&\  \hfil $##$\hfil\crcr
\mathstrut \crcr \noalign {\kern -\baselineskip } #1\crcr
\mathstrut \crcr \noalign {\kern -\baselineskip }}}\,}

\def\pmatrice#1{\left(\null\vcenter {\normalbaselines \m@th
\ialign {\hfil $##$\hfil &&\thinspace  \hfil $##$\hfil\crcr
\mathstrut \crcr \noalign {\kern -\baselineskip } #1\crcr
\mathstrut \crcr \noalign {\kern -\baselineskip }}}\right)}

\def\bmatrice#1{\left[\null\vcenter {\normalbaselines \m@th
\ialign {\hfil $##$\hfil &&\thinspace  \hfil $##$\hfil\crcr
\mathstrut \crcr \noalign {\kern -\baselineskip } #1\crcr
\mathstrut \crcr \noalign {\kern -\baselineskip }}}\right]}

\catcode`\@=12



\frenchspacing

\def\can{\mathop{\rm can}\nolimits}
\def\refCa{Ca72}
\def\refCF{CF69}
\def\refFo{Fo65}
\def\refFH{FH05}
\def\refGLZ{GLZ05}
\def\refGe{Ge87}
\def\refHo{Ho82}
\def\refMa{Ma15}
\def\refZ{Z80}

\rightline{2005/12/07 15:50}

\bigskip
\bigskip
\bigskip

\centerline{\bf A New Proof of the Garoufalidis-L\^e-Zeilberger}
\smallskip
\centerline{\bf Quantum MacMahon Master Theorem}
\medskip
\bigskip
\centerline{Dominique Foata and Guo-Niu Han}
\bigskip
\abstract{
We propose a new proof of the quantum version
of MacMahon's Master Theorem, established by Garoufalidis, L\^e and
Zeilberger. 
}
\medskip
\resume{
Nous proposons une nouvelle d\'emonstration de la version quantique
du Master Th\'eor\`eme de MacMahon,
\'etabli par Garoufalidis, L\^e et Zeilberger.
}

\bigskip
\centerline{\bf 1. Introduction} 

\medskip
MacMahon's Master Theorem is still regarded as a keystone in Combinatorial
Analysis. Numerous applications were already given by MacMahon himself
[\refMa, vol.~1, p.~97]. Several other proofs are due to Cartier [\refCa].
As the
Master Theorem is a special case of the multivariable  Lagrange Inversion
Formula [Ge87], as was shown by Hofbauer [\refHo], each
proof of the later formula yields a proof of this theorem. Its
first noncommutative version can be found in [\refFo],
which was later replaced into an appropriate algebraic
set-up [\refCF]. Recently Garoufalidis, L\^e and Zeilberger
have derived another noncommutative version [\refGLZ]
(called ``quantum Master Theorem") by using  difference
operator techniques developed by Zeilberger [\refZ]. In this
paper we propose a new proof of this quantum Master
Theorem.  As indicated in [\refFH], the formal parameter~$q$ 
in the Garoufalidis-L\^e-Zeilberger version plays no
role; we will then let $q=1$ in the present derivation
without loss of  generality.

\medskip

Let $r$ be a positive integer; the set
$\setA=\{1,2,\ldots,r\}$ is referred to as the underlying
{\it alphabet}.  A {\it biword} on $\setA$ is a $2\times n$
matrix $\alpha$ ($n\geq 0$), whose entries are in $\setA$,
the first (resp. second) row being called the {\it top word}
(resp. {\it bottom word}) of the biword $\alpha$. The
number $n$ is the length of $\alpha$; we write
$\len(\alpha)=n$.  Let $\cal B$ be the set of all biwords on
$\setA$. 
Each biword can also be viewed as a word of
{\it biletters} $\bw xa$  written vertically. The {\it product} of two
biwords is just the concatenation of them viewed as two
words of biletters. The biword of length $0$ is denoted by
$\bw{}{}=1$.

Let $\bboard Z$ be the ring of all integers. The set
$\setFBA={\bboard Z}\langle\!\langle {\cal B}\rangle\!\rangle$
of the formal sums $\sum_{\alpha} c(\alpha) \alpha$,
where $\alpha\in {\cal B}$ and $c(\alpha)\in\setZ$  for all
$\alpha\in{\cal B}$, together with the above multiplication,
the free addition and the free scalar product
is an algebra over $\setZ$, called the {\it free biword large
$\bboard Z$-algebra}.

\goodbreak
\medskip
The {\it right quantum algebra}~$\cal R$ is defined to be the
{\it associative} algebra, which is the quotient of ${\cal A}$
by the ideal~$\cal I$ generated by the commutation
relations
$$
\leqalignno{
 \textstyle{xy\choose aa} &=
   \textstyle{yx\choose aa},&\hbox{(R1)}\cr
 \noalign{\smallskip}
\textstyle{xy\choose ab} &=\textstyle
  {yx\choose ab} + {yx\choose ba}-{xy\choose ba}, &\hbox{(R2)} \cr
 }
$$
for all letters $a$, $b$, $x$, $y$ in~$\setA$ and $x\not= y$.
Notice that the associativity of the right quantum algebra
is set by {\it definition}. In fact, the
associativity is a {\it consequence} of the commutation
relations (R1) and (R2), as was proved in [\refFH]. All further
calculations in the paper will be made in the right quantum
algebra~$\cal R$, unless explicitly indicated. The elements
of~$\cal R$ will be called {\it expressions}. If ``can" is the
canonical homomorphism of~$\cal A$ onto~$\cal R$, we
identify $\can {\bboard B}$ with~$\bboard B$ and each
biword $x_1\;x_2\;\cdots \;x_m\choose x_{j_1}x_{j_2}\cdots
x_{j_m}$ occurring in an expression from~$\cal R$
with the product
${x_{i_1}\choose x_1}{x_{i_2}\choose x_2}\cdots
{x_{i_m}\choose x_m}$ in the quotient algebra~$\cal R$.

\medskip
For each word $w$ 
let $\overline w$ be its {\it non-decreasing}
rearrangement. The {\it Boson} is defined to be the infinite sum
$$
{\rm Bos}:=\sum_w  \pmatrice { \overline w\cr w\cr},
$$
where is the sum is over all words~$w$ from the free monoid
${\bboard A}^*$ generated by~$\bboard A$.
The {\it Fermion} is defined by
$$\Ferm
:=\sum_{J\subset \setAA}
(-1)^{|J|} \sum_{\sigma\in {\cal S}_J}
(-1)^{\inv \sigma}
\pmatrice  {\sigma (i_1)&\sigma (i_2)&\cdots&\sigma (i_l)\cr
i_1&i_2&\cdots&i_l\cr},
$$
where $J=\{i_1<i_2<\cdots <i_l\}$ and ${\cal S}_J$ is the
permutation group acting on the set~$J$. Both Boson and
Fermion will be regarded as elements of~$\cal R$ using the
above identification.
The {\it quantum Master Theorem} [\refGLZ] is stated next.

\medskip
\proclaim Theorem 1. The following identity
$$
\Ferm \times \Bos =1
$$
holds in the
right quantum algebra
${\cal R}$.

The proof of Theorem~1 is based on specific
techniques of Circuit Calculus (section 2) and determinantal
identities in the context of the right quantum algebra
(section~3). The end of the proof is made in
section 4.
 
\bigskip
\centerline{\bf 2. Real and imaginary expressions} 

\medskip
A {\it circuit} is a biword whose top word is a
rearrangement of its bottom word. Each formal sum
$E=\sum_\alpha c(\alpha)\alpha\in {\cal A}$ is said to be
{\it imaginary}  (resp. {\it real\/}), if $c(\alpha)=0$ for all
circuits $\alpha$ (resp. non-circuits $\alpha$).
If $E,E'\in {\cal A}$ and $E'\equiv E\!\!\!\pmod {\cal I}$,
then~$E'$ is imaginary (resp. real) if and only if~$E$ is
imaginary (resp. real), as easily verified by considering the
commutation relations (R1) and (R2). Accordingly, it makes
sense to say that an expression in~$\cal R$ is imaginary
(resp. real). Notice that $0$ is both real and imaginary. The
terminology is directly inspired from Complex Analysis.  The
following properties, although easy to verify, are essential
for the proof of the theorem.

\smallskip
(P1) Each expression $E\in{\cal R}$ can be decomposed in a unique way as
a sum of a real expression $\Re(E)$ and an imaginary expression $\Im(E)$:
$$E=\Re(E)+\Im(E).$$

\smallskip
(P2) The real part operator $\Re$ and the imaginary part operator $\Im$ are
linear. This means that for any two expresssions $E$ and
$E'$ we have
$$
\eqalign{  
\Re(E+E') & =\Re(E)+\Re(E'), \cr
\Im(E+E') & =\Im(E)+\Im(E'). \cr
}
$$

\smallskip
(P3) The real part operator $\Re$ and the imaginary part operator $\Im$ are
idempotent and orthogonal to each other. This means that for
every expresssion $E$ we have
$$
\Re(\Re(E))  =\Re(E),\
\Im(\Im(E))  =\Im(E)\
\hbox{and}\ 
\Re(\Im(E))  = \Im(\Re(E)= 0.
$$
\smallskip
(P4) If $E$, $E'$ are two nonzero expressions from~$\cal R$
and if $E$ is real, then
$E\times E'$ is real (resp. imaginary) if and only if $E'$ is
real (resp. imaginary).

\proclaim Lemma 2. If $E$, $E'$ are two nonzero expressions
from~$\cal R$ and if $E$ is real, then
$$\Re(E\times E')=E\times \Re(E').$$

{\it Proof}.\quad From the above properties we have
$$
\eqalign{  
 \Re(E\times E') &= \Re(E\times(\Re(E')+\Im(E'))) \cr
 &= \Re(E\times\Re(E')+E\times\Im(E')) \cr
 &= \Re(E\times\Re(E'))+\Re(E\times\Im(E')) \cr
 &= E\times\Re(E').\qed\cr
}
$$
\medskip

Now, define the {\it universe}  ``Univ" to be the sum of all
biwords whose top words are nondecreasing, that is,
$$
\Univ:=\sum_{u, w} {u\choose w},
$$
where $u$ (resp. $w$) runs over the set of all
{\it nondecreasing} words (resp. all words) and
where $u$ and $w$ are of the same length: $|u|=|w|$.
Both Boson and Fermion are sums of {\it circuits}. Moreover,
$\Re(\Univ)=\Bos$.
By Lemma~2 the quantum Master Theorem is equivalent to
the following theorem.

\proclaim Theorem 3. The following identity
$$
\Re(\Ferm \times \Univ) =1
$$
holds in the
right quantum algebra ${\cal R}$.


\section{3. Determinantal Calculus} 
Let $A=( a_{i,j} )_{1\leq i,j\leq r}$ be a square matrix
whose entries are expressions from~$\cal R$.
We define the {\it determinant} of~$A$ to be
$$
\det(A) = \
\sum_\sigma (-1)^{\inv \sigma}
a_{\sigma_1,1} a_{\sigma_2,2} \cdots a_{\sigma_r,r}.
$$
In this formula $\sigma$ runs over the permutation group
${\cal S}_{\setAA}$ of the set~$\setA$. The ordering of the
factors
$a_{\sigma_1,1}$, $a_{\sigma_2,2}$, $\ldots$, $a_{\sigma_r,r}$
{\it matters}, as the underlying algebra is noncommutative.
However, several classical properties of the determinant still
hold.

\proclaim Property 4 {\rm (Linearity)}.
When writing matrices as sequences of their $r$ columns
$(c_1,c_2,\ldots,c_r)$, we have
$$
\det(c_1, \ldots, c_i, \ldots, c_r)
+
\det(c_1, \ldots, c_i', \ldots, c_r)
=
\det(c_1, \ldots, c_i+c_i', \ldots, c_r).
$$

\proclaim Property 5 {\rm(Cofactor expansion)}.
Let $A=( a_{i,j} )_{1\leq i,j\leq r}$ be a matrix of
expressions and $A_{ij}$ be the matrix obtained from
$A$ by deleting the
$i$-th row and $j$-th column. Then
$$
\det(A)=\sum_{i=1}^r
(-1)^{r+i}\det(A_{ir}) a_{ir}.
$$

In the above identity we recognize the usual expansion of $\det(A)$ by
cofactors of the {\it rightmost} column.

\medskip
Let $a_i$, $b_i$, $c_i$ $(i=1,2,\ldots, r)$ and $x$, $y$ be
scalars and form the $r\times r$-matrix
$$
A=\pmatrix{\cdots&a_1+b_1{1\choose x}&c_1+b_1{1\choose y}&\ldots\cr
\cdots&a_2+b_2{2\choose x}&c_2+b_2{2\choose y}&\ldots\cr
\ddots&\vdots&\vdots&\ddots\cr
\cdots&a_r+b_r{r\choose x}&c_r+b_r{r\choose y}&\ldots\cr},
$$
where besides the two {\it consecutive} columns explicitly
displayed the other columns are arbitrary.  Let~$B$ denote the
matrix derived from~$A$ by transposing those two consecutive columns.

\proclaim Property 6 {\rm(Interchanging two
columns)}. Let $A$ and $B$ be the two matrices just
defined. Then
$\det A=-\det B$.

{\it Proof}.\quad
In the expansion of $\det A$ the sum of the two terms
$$
S_{ij}:=\pm \cdots (a_i+b_i\textstyle{i\choose x})(c_j+b_j{j\choose
y})\cdots
\ \mp \cdots (a_j+b_j\textstyle{j\choose x})(c_i+b_i{i\choose y})\cdots
$$
may be put in a one-to-one correspondence with the sum
$$
T_{ij}:=\pm \cdots (c_i+b_i\textstyle{i\choose y})(a_j+b_j{j\choose
x})\cdots
\ \mp \cdots (c_j+b_j{j\choose y})(a_i+b_i\textstyle{i\choose x})\cdots
$$
in the expansion of $\det B$. But
$$\displaylines{\quad S_{ij}=
\pm \cdots\Bigl(a_ic_j+a_ib_j\textstyle{j\choose y}+c_jb_i{i\choose
x}+b_ib_j{ij\choose
xy}\hfill\cr
\hfill{} -a_jc_i-a_jb_i\textstyle{i\choose y}-c_ib_j{j\choose
x}-b_ib_j{ji\choose
xy}\Bigr)\cdots\quad\cr
\quad T_{ij}=
\pm \cdots\Bigl(c_ia_j+c_ib_j\textstyle{j\choose x}+b_ia_j{i\choose
y}+b_ib_j{ij\choose
yx}\hfill\cr
\hfill{} -c_ja_i-c_jb_i\textstyle{i\choose x}-b_ja_i{j\choose
y}-b_ib_j{ji\choose
yx}\Bigr)\cdots.\quad\cr}
$$
Because of the identity
${ij\choose xy}-{ji\choose xy}=-({ij\choose yx}-{ji\choose yx})$
and the fact that the commutation relations can be made at {\it any}
position
within each biword (using the associativity property of the
right quantum algebra) we conclude that
$S_{ij}=-T_{ij}$.\qed
\medskip

It is worth noticing that Property 6 is only stated for very
special matrices. In general, the column interchanging
property does {\it not\/} hold for arbitrary matrices with
entries in~$\cal R$.
\medskip
Consider the $r\times r$-matrix
$\matB^{(r)}=\bigl(\, \bw{i}{j}\, \bigr)_{1\leq i,j\leq r}$.
Notice that every entry in $\matB^{(r)}$ is a biletter.
Define the {\it fermion-matrix} to be the matrix
$\matF^{(r)}=\matI-\matB^{(r)}$.
Using Property~4 (linearity) we have
$$
\Ferm = \det(\matF^{(r)})
=\det\pmatrix{
1-\bw 11 & -\bw 12 &  \ldots & -\bw 1r \cr
-\bw 21 & 1-\bw 22 &  \ldots & -\bw 2r \cr
\vdots & \vdots & \ddots & \vdots \cr
-\bw r1 & -\bw r2 &  \ldots & 1-\bw rr \cr
}. 
$$

Let $i\leq r-1$ and replace
the rightmost column in $\matF^{(r)}$
by the $i$-th column of $\matF^{(r)}$.
Let $\matF_i$ 
be the resulting matrix,
so that $\matF_i$ has two identical columns.

\proclaim Lemma 7.
We have: $\det(\matF_i)=0$.
 
{\it Proof}.\quad
Permute columns $i$ and $i+1$, then columns $i+1$ and $i+2$,
$\cdots$, finally columns $r-2$ and $r-1$,
We obtain a matrix $A$  whose rightmost two columns are identical.
Property~6 implies that $\det(A)=0$ and also
$\det(\matF_i)=\pm\det(A)=0$.\qed

\goodbreak
\bigskip
\centerline{\bf 4. The proof} 

\medskip
First, define
$$
\eqalign{  
S_i=S_i^{(r)}&:=\bw i1 +\bw i2 + \cdots +\bw ir;\cr
K_i=K_i^{(r)}&:= {1\over 1-S_i}=\sum_{n\geq 0} S_i^n.\cr
}$$

\proclaim Lemma 8.
The universe defined in Section $2$ is equal to
$$
\Univ=K_1 K_2 \cdots K_r.
$$

{\it Proof}.\quad
By definition of ``$\Univ$" we have
$$
\eqalign{  
\Univ&=\sum_{u,w} \bw {u}w=
\sum_{w_1,w_2,\cdots,w_r} \bw {11\cdots1}{w_1}
\bw {22\cdots2}{w_2}
\cdots
\bw {rr\cdots r}{w_r}\cr
&=
\sum_{w_1} \bw {11\cdots1}{w_1}
\sum_{w_2}\bw {22\cdots2}{w_2}
\cdots
\sum_{w_r}\bw {rr\cdots r}{w_r}\cr
&=K_1K_2\cdots K_r. \qed\cr
}$$

\proclaim Lemma 9.
We have: $S_iS_j=S_jS_i$ and $K_iK_j=K_jK_i$.

{\it Proof}.
Grouping the biwords by pairs if necessary we have:
$$
\eqalign{  
S_iS_j&=\sum_{a<b} (\bw{ij}{ab}+\bw {ij}{ba}) + \sum_a \bw{ij}{aa} \cr
&=\sum_{a<b} (\bw{ji}{ab}+\bw{ji}{ba})+\sum_a \bw{ji}{aa}\cr
&=S_jS_i.\qed\cr
}$$

{\it Proof of Theorem $3$}.
Let $\matM$ be the matrix obtained from the
fermion-matrix $\matF^{(r)}$ by adding
all the leftmost $r-1$ columns to the rightmost column:
$$
\matM= \pmatrix{
1-\bw 11 & -\bw 12 &  \ldots & 1-S_1 \cr
-\bw 21 & 1-\bw 22 &  \ldots & 1-S_2 \cr
\vdots & \vdots & \ddots & \vdots \cr
-\bw r1 & -\bw r2 &  \ldots & 1-S_r \cr
}  
$$
By Property 4 and Lemma 7: $\det(\matM)=\det(\matF^{(r)})=\Ferm$.
Let
$$
E^{(r)}:=\Ferm\times \Univ =\det(\matM)\times \Univ
$$
By Lemmas 8 and 9
$$
\eqalign{  
E^{(r)}
&= \det\pmatrix{
1-\bw 11 & -\bw 12 &  \ldots & (1-S_1)\times \Univ \cr
-\bw 21 & 1-\bw 22 &  \ldots & (1-S_2)\times \Univ \cr
\vdots & \vdots & \ddots & \vdots \cr
-\bw r1 & -\bw r2 &  \ldots & (1-S_r)\times \Univ \cr
}  \cr
&= \det\pmatrix{
1-\bw 11 & -\bw 12 &  \ldots & K_2K_3\cdots K_r\cr
-\bw 21 & 1-\bw 22 &  \ldots & K_1K_3\cdots K_r \cr
\vdots & \vdots & \ddots & \vdots \cr
-\bw r1 & -\bw r2 &  \ldots & K_1K_2\cdots K_{r-1}\cr
}.  \cr
}$$
Now applying Property 5 to the above determinant yields
$$
\eqalign{
E^{(r)} = & (-1)^{r+1} \det(\matM_{1r}) K_2K_3\cdots K_r \cr
      & (-1)^{r+2} \det(\matM_{2r})  K_1K_3\cdots K_r +\cdots \cr
      & + \det(\matM_{rr}) K_1K_2\cdots K_{r-1}, \cr
}
$$
where, as in Property 5, $\matM_{ij}$ denotes the minor at
position $(i,j)$. For each biword occurring in
$F_1:=\det(\matM_{1r}) K_2K_3\cdots K_r$ there is a letter
$1$ in the bottom, but none in the top. This means that
$F_1$ does not contain any circuit. Thus $\Re (F_1) =0$.
This argument is also valid for the other terms in the above
summation of $E^{(r)}$, except for the last term.
When taking the real part of the $E^{(r)}$ we obtain:
$$
\eqalign{  
\Re(E^{(r)})
&=\Re (\det(\matM_{rr}) K_1K_2\cdots K_{r-1}) \cr
&=\Re (\det(\matF^{(r-1)}) K_1^{(r-1)}K_2^{(r-1)}\cdots K_{r-1}^{(r-1)}) \cr
&=\Re (E^{(r-1)}).\cr
}$$
Then, by iteration
$$\Re(E^{(r)})= \Re(E^{(r-1)})=
\cdots  =\Re(E^{(1)})
= \Re((1-\bw 11) K_1^{(1)})
=1. \qed
$$

\vfill\eject\vglue 15mm
\bigskip

{
\eightpoint

\centerline{\bf References} 
\bigskip

\livre \refCa|P. Cartier|
La s\'erie g\'en\'eratrice exponentielle, Applications probabilistes
et alg\'ebriques|
Publ. I.R.M.A., Universit\'e Louis Pasteur, Strasbourg,
{\oldstyle 1972}, 68 pages|
\smallskip

\livre \refCF|P. Cartier, D. Foata|Probl\`emes combinatoires de
permutations et r\'earran\-ge\-ments|Berlin, Springer-Verlag,
{\oldstyle 1969} ({\sl Lecture Notes in Math., {\bf 85}})|
\smallskip

\article \refFo|D. Foata|Etude alg\'ebrique de certains probl\`emes
d'Analyse Combinatoire et du Calcul des Probabilit\'es|Publ.
Inst. Statist. Univ. Paris|14|1965|81--241|
\smallskip

\divers \refFH|D. Foata, G.-N. Han|{\sl A basis for the right quantum
algebra
and the ``$1=q$" principle}, preprint, {\oldstyle 2005}|
\smallskip

\divers \refGLZ|S. Garoufalidis, T. Tq L\^e, D. Zeilberger|
{\sl The Quantum MacMahon Master Theorem}, to be appear
in {\sl Proc. Natl. Acad. Sci.},
{\oldstyle 2005}|
\smallskip

\article \refGe|I. M. Gessel|A Combinatorial Proof of the Multivariable
Lagrange Inversion Formula|J. Comb. Theory, Ser. A|45|1987|178--195|
\smallskip

\divers \refHo|J. Hofbauer|Lagrange Inversion, {\sl S\'eminaire Lotharingien
de Combinatoire}, {\bf B06a}, {\oldstyle 1982}, 28 pages,
({\tt http://www.mat.univie.ac.at/$\!\sim$slc})|
\smallskip
 
\livre \refMa|P. A. MacMahon|Combinatory Analysis {\rm vol.
I, II}|Cambridge Univ. Press, {\oldstyle 1915}. Reprinted by
Chelsea, New York, {\oldstyle 1960}|
\smallskip

\article \refZ|D. Zeilberger|The algebra of linear partial difference
operators and its applications|SIAM J. Math. Anal.|11, no. 6|1980|919--932|
\smallskip


} 

\bigskip\bigskip
\hbox{\vtop{\halign{#\hfil\cr
Dominique Foata \cr
Institut Lothaire\cr
1, rue Murner\cr
F-67000 Strasbourg, France\cr
\noalign{\smallskip}
{\tt foata@math.u-strasbg.fr}\cr}}
\qquad
\vtop{\halign{#\hfil\cr
Guo-Niu Han\cr
I.R.M.A. UMR 7501\cr
Universit\'e Louis Pasteur et CNRS\cr
7, rue Ren\'e-Descartes\cr
F-67084 Strasbourg, France\cr
\noalign{\smallskip}
{\tt guoniu@math.u-strasbg.fr}\cr}}
}

\vfill\eject

\bye